\documentclass[12pt]{article}
\usepackage[vmargin=1in,hmargin=1in]{geometry}

\usepackage{amssymb, cite}
\usepackage{graphicx}

\newcommand{\be}{\begin{equation}}
\newcommand{\ee}{\end{equation}}
\newcommand{\bqn}{\begin{eqnarray}}

\newcommand{\eqn}{\end{eqnarray}}

\newcommand{\bd}{\begin{description}}
\newcommand{\ed}{\end{description}}
\newtheorem{theorem}{Theorem}[section]
\newtheorem{lemma}{Lemma}[section]

\newtheorem{question}{Question}[section]

\newtheorem{stat}{}[section]

\def\bs{\begin{stat}}
\def\es{\end{stat}}
\def\ben{\begin{enumerate}}
\def\een{\end{enumerate}}

\def\bp{\noindent{\bf Proof}  \ }
\newcommand{\ep}{\hfill $\square$}

\makeatletter
\@addtoreset{equation}{section}

\makeatother

\begin{document}

\begin{center}
European Journal of Combinatorics, 34 (2013) 764 - 769
\end{center}

\begin{center}
{\large {\bf ON  CONVEX POLYTOPES  IN $\mathbb{R}^d$ 
\\[0.7ex]
CONTAINING AND AVOIDING
ZERO}}
\\[4ex]
{\large {\bf Alexander Kelmans}}
\\[2ex]
{\bf University of Puerto Rico, San Juan, Puerto Rico, United States}
\\[0.5ex]
{\bf Rutgers University, New Brunswick, New Jersey, United States}
\\[4ex]
{\large {\bf Anatoly Rubinov}}
\\[2ex]
{\bf Institute for Information Transmission Problems, Moscow, Russia}
\end{center}

\footnotetext
{{\em E-mail address}: kelmans.alexander@gmail.com 
(Alexander Kelmans)}

\begin{abstract}
The goal of this paper is to establish certain inequalities between the numbers of convex polytopes in 
$\mathbb{R}^d$ ``containing'' and ``avoiding'' zero provided that their vertex sets are subsets of a given finite set $S$ of points in $\mathbb{R}^d$. 
This paper is motivated by a question about these quantities raised by E. Boros and V. Gurvich in 2002.
\\[1ex]
\end{abstract}

\section{Introduction}

\indent

The notions and facts used but not described here can be found in \cite{S}.

Let $S$ be a finite set of points in $\mathbb{R}^d$, $X \subseteq S$,  and 
$z \in \mathbb{R}^d \setminus S$.
\\[1ex]
A set $X$ is called a {\em $z$-containing set} (a  {\em $z$-avoiding set}) if $z$ is in the interior of the convex hull of $X$ 
(respectively, $z$ is not in the interior of the convex hull of $X$). 
\\[1ex]
A $z$-containing set $X$ is {\em minimal}  if $X$ has no proper $z$-containing subset.
\\[1ex]
A $z$-avoiding set $X$ is {\em maximal in}  $S$ if $S$ has no $z$-avoiding subset containing $X$ properly.
\\[1ex]
Let ${\cal C}(S)$ and ${\cal A}(S)$  denote the sets of minimal 
$z$-containing and maximal  $z$-avoiding subsets of $S$, respectively.
\\[1ex]
\indent
In 2002 E. Boros and V. Gurvich raised the following interesting question.
\begin{question}
\label{question} Suppose that $S$ is a finite set of points in $\mathbb{R}^d$,
 $z \in \mathbb{R}^d \setminus S$,  and $S$  is a $z$-containing set. Is it true that  
$|{\cal A}(S)| \le 2 d~|{\cal C}(S)|$?
\end{question}
%

Questions of this type arise naturally in the algorithmic theory of a so-called {\em efficient} enumeration of different type of geometric or combinatorial objects 
(see, for example, \cite{Rtcr, IIS}).
\\[1ex]
\indent
 For an affine subspace $F$ of $\mathbb{R}^d$, let $dim(F)$ denote the affine dimension of $F$ 
and  $R(X)$ denote the minimal affine subspace in $\mathbb{R}^d$ containing $X$.
\\[1ex]
\indent
A finite set $S$ of points in $\mathbb{R}^d$ is said to be {\em in a general position} if for every $X \subseteq S$, 
\\[0.3ex]
$|X| - 1 \le d$ $\Rightarrow$ $dim (R(X)) = |X| - 1$.
\\[1ex]
\indent
We say that {\em $z$  is in a general position with respect to $S$} if
$dim(R(X)) < d$ $\Rightarrow$ $z \not \in R(X)$ for every $X \subseteq S$.
\\[1ex]
\indent
One of our main results is the following theorem.
\begin{theorem} 
\label{!}
Let $S$ be a finite set of points in the 
$d$-dimensional space $\mathbb{R}^d$, $z \in \mathbb{R}^d \setminus S$.
Suppose that $z$  is in a general position with respect to $S$.
Then  $|{\cal A}(S)| \le d~|{\cal C}(S)| + 1$.
\end{theorem}

This theorem was announced in \cite{CKK} and its proof was presented at the RUTCOR seminar directed by E. Boros and V. Gurvich in June 2002.
\\[1ex]
\indent
Later L. Khachian gave a construction providing for every 
$d \ge 4$ a counterexample $(S,z)$ in  $\mathbb{R}^d$ to the inequality in Question \ref{question} such that $S$ is not in a general position (see \cite{Rtcr}).
From Theorem \ref{!} it follows that in all these counterexamples $z$ is not in a general position with respect to $S$.
\\[1ex]
\indent
In \cite{CKK} it is shown that if $S$ is a finite set of points in the plane 
$\mathbb{R}^2$, $z \in \mathbb{R}^2 \setminus S$, and $S$ is a $z$-containing set, then 
$|{\cal A}(S)| \le 3|{\cal C}(S)| +1$, and so 
the inequality in Question \ref{question} is true for the plane.
\\[1ex]
\indent
In this paper we also give some strengthenings of Theorem  \ref{!}.


\section{Some notions, notation, and auxiliary facts}
\label{notions}

\indent

Given a convex polytope $P$ in $\mathbb{R}^d$, a  face $F$ of a polytope $P$ is called 
a  {\em facet} of $P$ if $dim (R(F)) = d - 1$.
Obviously, $R(F)$ is a hyperplane; we call $R(F)$ 
a {\em facet hyperplane of $P$}.
Let, as above, $S$ be a finite set of points in $\mathbb{R}^d$, $X \subseteq S$,  and 
$z \in \mathbb{R}^d \setminus S$. Let $s \in S$.
\\[1ex]
\indent
We will use the following notation:

\begin{itemize}
\item $conv(X)$ is the convex hull of $X$,

\item if $P$ is a convex polytope, then $V(P)$ is the set of vertices of $P$ and $v(P) = |V(P)|$, 

\item as above, ${\cal C}(S)$ is the set of minimal 
$z$-containing  subsets of $S$; 
also ${\cal C}_s(S)$ is the set of members of ${\cal C}(S)$ containing $s$,

\item as above, ${\cal A}(S)$  is the set of  maximal  
$z$-avoiding subsets of $S$,

\item ${\cal A}^s(S)$ is the set of subsets $A$ in $S$ such that $A$ is maximal $z$ avoiding in $S$ and 
$A \setminus s$ is maximal $z$-avoiding in 
$S \setminus s$, and so $A \in {\cal A}(S)$ and 
$A \setminus s \in {\cal A}(S \setminus s)$,

\item ${\cal A}_s(S) = {\cal A}(S) \setminus {\cal A}^s(S)
 = \{X \in {\cal A}(S): s \in X~and~X \setminus s \not \in {\cal A}(S \setminus s\}$,

\item ${\cal S}mpl(S)$ is the set of simplexes $C$ such that
$z$ is an interior point of $C$ and $V(C) \subseteq S$ 
and
${\cal S}mpl_s(S)$ is the set of simplexes $C$ in 
${\cal S}mpl(S)$ such that $s \in  V(C)$, 

\item ${\cal C}onv(S) = \{conv(X): X \in {\cal C}(S)\}$
and ${\cal C}onv_s(S)$ is the set of members of 
${\cal C}onv(S)$ containing $s$, 

\item  ${\cal H}(S)$ is the set of hyperplanes $H$ such that $H$ is a facet hyperplane of a simplex in  ${\cal S}mpl(S)$ and ${\cal H}_s(S)$ is the set of hyperplanes in ${\cal H}(S)$ containing $s$, and

\item ${\cal F}(S)$ is the set of subsets $T$ of $S$ such that $|T| = d$ and $conv(T)$  is a face of a simplex in  ${\cal S}mpl(S)$ and ${\cal F}_s(S)$ is the set of subsets of $S$ in ${\cal F}(S)$ containing $s$.

\end{itemize}

Obviously, 
we have the following.
\begin{lemma}
\label{observation}
Let $S$ be a finite set of points in $\mathbb{R}^d$.
Then
\\[1ex]
$(a1)$ $|{\cal C}onv(S)| = |{\cal C}(S)|$,
\\[1ex]
$(a2)$
 $|{\cal H}(S)| \le |{\cal F}(S)| \le (d+1)|{\cal S}mpl(S)|$, and 
 $|{\cal C}onv_s(S)| = |{\cal C}_s(S)|$,
\\[1ex]
$(a3)$
$|{\cal H}_s(S)| \le |{\cal F}_s(S)| \le d~|{\cal S}mpl_s(S)|$, 
and
\\[1ex]
$(a4)$
$|{\cal A}^s(S)| = |{\cal A}(S\setminus s)|$, and so 
$|{\cal A}(S)| -  |{\cal A}(S\setminus  s)| = |{\cal A}_s(S)|$.\end{lemma}

We recall that {\em $z$  is in a general position with respect to $S$} if
$dim(R(X)) < d$ $\Rightarrow$ $z \not \in R(X)$ for every $X \subseteq S$.
\\[1ex]
\indent
We will use the following well known and intuitively obvious fact.
\begin{lemma}
\label{SimplexInPcontainingZero}
Let $P$ be a convex polytope in $\mathbb{R}^d$ and 
$z$ a point in the interior of $P$. Suppose that
$z$ is in general position with respect to $V(P)$.
Then there exists $X \subseteq V(P)$ such that
$conv(X)$ is a simplex of dimension $d$ and $z$ is in the interior of  $conv(X)$.
\end{lemma}

\bp 
We prove our claim by induction on the dimension $d$.
The claim is obviously true for $d = 1$.
We assume that the claim is true for $d = n-1$ and will prove that the claim is also true for 
$d = n$, where $n \ge 2$. Thus, $P$ is a convex polytope in $\mathbb{R}^n$.
Since $z$ is in a general position with respect to $V(P)$, clearly $z \not \in V(P)$. Let $p \in V(P)$ and $L$ the line containing $p$ and $z$. Then there exists the point $t$ in $L \cap P$ such that the closed interval
$pLt$ contains $z$ as an interior point and $t$ is not an interior point of $P$.
In a plain language, $pL$ is the ray going from point $p$ through point $z$ and $t$ is the first point of the ray which is  not the interior point of $P$. Then $t$ belongs to a face $F$ of $P$ which is a convex polytope of dimension at most $n - 1$. Since $z$ is in a general position with respect to $V(P)$, the dimension of  $F$ is  $n-1$ and $t$ is in a general position with respect to $V(F)$. By the induction hypothesis, there exists 
$T \subseteq V(F)$ such that  $conv(T)$ is a simplex of dimension $n - 1$ and $t$ is in the interior of  
$conv(T)$. Put $X = T\cup p$. Then $X \subseteq V(P)$, 
$conv(X)$ of dimension $n$ is a simplex and $z$ is in the interior of  $conv(X)$.
\ep
\\[1.5ex]
\indent
Lemma \ref{SimplexInPcontainingZero} also follows from the Caratheodory Theorem (see \cite{S}) and Lemma \ref{z-in-conv(X)} below.
\\[1.5ex]
\indent
Given $X, Y \subset \mathbb{R}^d$ and a hyperplane $H$, we say that $H$ {\em separates $X$ and $Y$} (or {\em separates $X$ from $Y$}) if
$X \setminus H$ and $Y \setminus H$ belong to different 
half-spaces of $\mathbb{R}^d \setminus H$.
\\[1.5ex]
\indent
It is easy to see the following.
\begin{lemma}
\label{(A,z)-separating-hyperplane}
Let $A$ and $A'$ be  $z$-avoiding subsets of $S$. Then
\\[1ex]
$(a1)$ there exists a hyperplane separating $A$ from $z$,
\\[1ex]
$(a2)$
if there exists a hyperplane separating both $A$ and $A'$ from $z$, then $A \cup A'$ is also a $z$-avoiding subsets of $S$, and therefore  
\\[1ex]
$(a3)$ if $A$ and $A'$ are maximal $z$-avoiding subsets of $S$ and there exists a hyperplane separating both $A$ and $A'$ from $z$, then $A= A'$.

\end{lemma}

We also need the following simple facts.
\begin{lemma}
\label{z-in-conv(X)}
Let $z$  be in a general position with respect to 
$S$ and $X \subseteq S$.
Then the following are equivalent:
\\[0.7ex]
$(a1)$ $z \in conv(X)$ and
\\[0.7ex]
$(a2)$ $z$ is in the interior of $conv(X)$.
\end{lemma}

\bp 
Obviously, $(a2)$ implies $(a1)$.
Suppose, on the contrary,  $(a1)$ does not imply $(a2)$.
Then $z$ belongs to  $conv(X')$ for some 
$X' \subseteq X$ with $dim (R(X')) < d$. 
Therefore  $z$  is not in a general position with 
respect to $S$, a contradiction.
\ep

\begin{lemma}
\label{Min-z-containingSet}
Let $z$  be in a general position with respect to 
$S$. If $C$ is a minimal $z$-containing subset of $S$,
then $conv(C)$ is a simplex, and so
${\cal S}mpl(S) = {\cal C}onv(S)$. 
\end{lemma}

\bp (uses Lemmas 
\ref{SimplexInPcontainingZero} and \ref{z-in-conv(X)}).
Since $C$ is $z$-containing, there exists $X \subseteq C$ such that $z \in conv(X)$ and $conv(X)$ is a simplex.
By Lemma \ref{z-in-conv(X)}, $z$ is in the interior of $conv(X)$.
Since $C$ is minimal $z$-containing, clearly 
$C = X$.
Now by Lemma 
\ref{SimplexInPcontainingZero}, $conv(C)$ is a simplex.
\ep

\begin{lemma}
\label{sz-ray}
Let $z$  be in a general position with respect to $S$. Suppose that $X \subset S$ and $s \in S \setminus X$.
Let $L$ be the line in $\mathbb{R}^d$ containing $s$ and $z$.
If $dim (R(X)) \le d-2$, then  $L \cap R(X) = \emptyset $.
\end{lemma}

\bp Suppose, on the contrary, $L \cap R(X)\ne \emptyset $.
Then $z \in R(X \cup s)$ and $dim (R(X \cup s)) < d$. 
Then $z$ is not in a general position with respect to $S$, 
a contradiction.
\ep

\begin{lemma}
\label{sLz-A,A'}
Let $S$ be a $z$-containing set and $z$   in a general position with respect to $S$.
Let $A$ be a maximal $z$-avoiding set in $S$ and
$s \in S \setminus A$.
Then $A$ has a subset $T$ such that
\\[1ex]
$(a1)$
$|T| = d$, 
\\[1ex]
$(a2)$
$T$ belongs to a facet of $conv(A)$, 
\\[1ex]
$(a3)$  
$conv(T \cup s)$ is a simplex of dimension $d$ containing $z$ as an interior point, and 
\\[1ex]
$(a4)$
$R(T)$ is a facet hyperplane of $A$ separating
$A$ from  $z$.
\end{lemma}


\bp (uses Lemmas \ref{SimplexInPcontainingZero} and  \ref{sz-ray}).
Let $L$ be the line in $\mathbb{R}^d$ containing $s$ and $z$. 
Since $s \not \in A$ and $A$ is a maximal $z$-avoiding subset in $S$, clearly $z$ is  in the interior of 
$conv(A \cup s)$.

Let $I = L \cap conv(A \cup s)$. Since $conv(A \cup s)$ is a convex set, clearly $I$ is a  line segment $sLr$, where $r$ 
is not an interior point of $conv(A \cup s)$.
We claim that $r \in conv(A)$. Indeed, if not, then $r$ belongs to a facet $R$ of $conv(A \cup s)$ containing $s$.
Since $sLr$ is a convex set and $s, r \in R$, we have:
$z \in sLr \subseteq R$, and so $z \in R$.
It follows that $z$ is not an interior point of   $conv(A \cup s)$, a contradiction.
Thus,   $conv(A) \cap L\ne \emptyset $.

Since $s \not \in A$, 
there exists a unique point $t$ in $L$ such that 
the closed interval $sLt$ in $L$ with the end-points $s$ and $t$ has the properties: $z$ is in the interior of $sLt$  
 and  $conv(A) \cap sLt = t$. 
 In a plane language, $sL$ is the ray going from point $s$ through point $z$ and $t$ is the first point of the ray belonging to $conv(A)$. 
By Lemma \ref{sz-ray}, $t$ belongs to the interior of a 
facet $F$ of $conv(A)$.
Hence, by Lemma 
\ref{SimplexInPcontainingZero}, $F$ has a set $T$ of 
$d$ vertices such that $T \subset A$ and 
$conv (T)$ is a simplex of dimension $d - 1$ containing 
$t$ as an interior point.
Then $conv (T \cup s)$ is a simplex of dimension $d$ containing 
$z$ as an interior point, and so $t$ is an interior point of $conv(T)$. 
\ep

\section{Main results}

\indent

First we will  prove a weaker version of our main result
to demonstrate the key idea concerning the relation between the minimal $z$ containing and maximal
$z$-avoiding subsets of $S$.
\begin{theorem}
\label{w}
Suppose that $S$ is a $z$-containing set and $z$ is in a general position with respect 
\\[0.5ex]
to $S$.
Then $|{\cal A}(S)| \le (d +1) |{\cal C}(S)|$.
\end{theorem}

\bp (uses Lemmas
\ref{(A,z)-separating-hyperplane}$(a4)$, 
\ref{Min-z-containingSet}, and
 \ref{sLz-A,A'}).
Let $A$ be a maximal $z$-avoiding subset of $S$.
Since $S$ is a $z$-containing set, there exists $s \in S \setminus A$.
Since $A$ is a maximal $z$-avoiding subset of $S$,
$A \cup s$ is not a $z$-avoiding subset of $S$. Hence 
$A \cup s$ is a $z$-containing set, and so  $z$ is  in the interior of $A \cup s$.
By Lemma \ref{sLz-A,A'}, 
$A$ has a subset $T$ such that
$|T| = d$, $T$ belongs to a face of $conv(A)$,  
and 
$conv(T \cup s)$ is a simplex containing $z$ as an interior point.
\\[0.7ex]
\indent
Let, as above,  ${\cal H}(S)$ denote the set of hyperplanes $H$ such that $H$ is a facet hyperplane of a simplex in  ${\cal S}mpl(S)$ and ${\cal F}(S)$ denote the set of subsets $T$ of $S$ such that $|T| = d$ and $conv(T)$  is a facet of a simplex in  ${\cal S}mpl(S)$.
Obviously, $|{\cal H}(S)| \le |{\cal F}(S)| = (d+1)|{\cal S}mpl(S)|$.
By  Lemmas \ref{(A,z)-separating-hyperplane}$(a4)$ and \ref{sLz-A,A'}, $|{\cal A}(S)| \le |{\cal H}(S)|$.
Clearly,
$|{\cal C}onv(S)| = |{\cal C}(S)| $ and, by Lemma 
\ref{Min-z-containingSet},
${\cal S}mpl(S) = {\cal C}onv(S)$. 
Thus, $|{\cal A}(S)| \le  (d+1) |{\cal C}(S)| $.
\ep
\\[1.5ex]
\indent
One of the referees informed us that Theorem \ref{w}
was formulated in terms of minimal infeasible subsystems  and proved in a different way in \cite{IIS}.
\\[1.5ex]
\indent
A hyperplane $H$ in ${\cal H}(S)$ is said to be 
{\em essential}
if $H$ is a facet hyperplane of a maximal $z$-avoiding subset $A$ in $S$ separating $A$ from $z$, and
{\em non-essential}, otherwise.
Let ${\cal H}^e(S)$ and ${\cal H}^e_s(S)$ denote the sets of essential hyperplanes in ${\cal H}(S)$ and ${\cal H}_s(S)$, respectively.
\begin{lemma}
\label{lemma!}
Let $S$ be a finite set of points in $\mathbb{R}^d$, $z \in \mathbb{R}^d \setminus S$, and $s \in S$.
Suppose that $S$ is a $z$-containing set and $z$  is in a general position with respect to $S$.
Then 
\\[1ex]
$|{\cal A}(S)| -  |{\cal A}(S\setminus  s)| = 
|{\cal A}_s(S)| \le |{\cal H}^e_s(S)| \le |{\cal H}_s(S)| \le |{\cal F}_s(S)|\le 
 d~|{\cal S}mpl_s(S)| =d~|{\cal C}_s(S)|$.
\end{lemma}

\bp (uses  Lemmas   
\ref{(A,z)-separating-hyperplane}$(a4)$,
\ref{Min-z-containingSet},
 and \ref{sLz-A,A'}).
We prove that  $|{\cal A}_s(S)| \le |{\cal H}^e_s(S)|$.
Let $A \in {\cal A}_s(S)$. Then $s \in A $ and 
$A' = A \setminus s$ is a $z$-avoiding but not maximal 
$z$-avoiding  set in $S' = S - s$. Therefore there exists 
$s' \in S' \setminus A'$ such that $A' \cup s'$ is also a 
$z$-avoiding set. 
Obviously, $A' \cup \{s,s'\}$ is a $z$-containing set in $S$.
By Lemma \ref{sLz-A,A'}, 
$A$ has a subset $T$ such that
$|T| = d$, $T$ belongs to a face of $conv(A)$,  
and 
$conv(T \cup s')$ is a simplex containing $z$ as an interior point. Since $A' \cup s'$ is  a $z$-avoiding set, clearly
$s \in T$.
Now by Lemmas   \ref{(A,z)-separating-hyperplane}$(a4)$ and \ref{sLz-A,A'}, $|{\cal A}_s(S)| \le |{\cal H}^e_s(S)|$.
By Lemma \ref{Min-z-containingSet},
${\cal S}mpl_s(S) = {\cal C}onv_s(S)$, and clearly,
$|{\cal C}onv_s(S)| = |{\cal C}_s(S)| $.
All the other inequalities in our claim are obvious.
\ep.
\\[1.5ex]
\indent
Now we are ready to prove the following strengthening of  
Theorem \ref{w} which is also an extension of 
Theorem \ref{!}.
\begin{theorem} 
\label{m}
Let $S$ be a finite set of points in the 
$d$-dimensional space $\mathbb{R}^d$.
Suppose that $z$  is in a general position with respect to $S$ {\em (and so $z \in \mathbb{R}^d \setminus S$)}.
Then 
\\[1ex]
$(a1)$ if  either 
$S$ is $z$-avoiding  or $S$ is $z$-containing and $|S| = d+1$ {\em (and so $conv(S)$ 
is a $d$-dimensional simplex)}, then 
$|{\cal A}(S)| = d~|{\cal C}(S)| + 1$,  
\\[1ex]
$(a2)$ $|{\cal A}(S)| \le d~|{\cal C}(S)| + 1$, 
\\[1ex]
$(a3)$  if $S$ is $z$-containing and $|S| = d + 2$, then
$|{\cal A}(S)| = d~|{\cal C}(S)| - d +1$, and
\\[1ex]
$(a4)$ if $S$ is $z$-containing and $|S| \ge d + 3$,
then $|{\cal A}(S)| \le d~|{\cal C}(S)| - d$.
\end{theorem} 
%

\bp(uses Lemmas 
\ref{Min-z-containingSet}  
and 
\ref{lemma!}).
\\[1.5ex]
\indent
{\bf (p1)} First we prove $(a1)$. If $S$ is $z$-avoiding, then 
$|{\cal A}(S)| = 1$ and   $|{\cal C}(S)| = 0$, and so 
$|{\cal A}(S)| = d~|{\cal C}(S)| + 1$.
If  $S$ is $z$-containing  and $|S| = d+1$, then
$conv(S)$ is a $d$-dimensional simplex, $z$ is in the interior of $conv(S)$,  $|{\cal A}(S)| = d+1$, and $|{\cal C}(S)| = 1$, and therefore 
$|{\cal A}(S)| = d~|{\cal C}(S)| + 1$.
\\[1.5ex]
\indent
{\bf (p2)}
We prove $(a2)$.
Our claim is obviously true if $S$ is a $z$-avoiding set.
Therefore we assume that $S$ is a $z$-containing set.
We prove our claim by induction on $|S|$.
By Lemma \ref{Min-z-containingSet}, $|S| \ge d+1$.
If $|S| = d+1$, then $conv(S)$ is a simplex, and our claim is obviously true. Thus, we assume that our claim is true for every $z$-containing set $S$ with $|S| = k \ge d+1$ and will prove that the claim is also true if $|S| = k + 1$.
Since $S$ is $z$-containing, by Lemma \ref{Min-z-containingSet}, there exists $X \subseteq S$ such that $conv(X)$ is a simplex and $z$ is the interior point of 
$conv(X)$. Since $|X| = d+1 < |S|$, there exists 
$s \in S \setminus X$. Obviously, $S' = S \setminus s$ is a $z$-containing set. Since $k = |S'| < |S| = k+1$, by the induction hypothesis, our claim is true for $S' = S \setminus s$, i.e.  $|{\cal A}(S')| \le d~|{\cal C}(S')| + 1$.
By Lemma \ref{lemma!}, 
$|{\cal A}(S)| -  |{\cal A}(S')| = 
|{\cal A}_s(S)| \le d~|{\cal C}_s(S)| $.
Now since $|{\cal C}(S)| =  |{\cal C}(S')| + |{\cal C}_s(S)|$, 
our inductive step follows.
\\[1.5ex]
\indent
{\bf (p3)} We prove $(a3)$. 
Since $S$ is $z$-containing, there exists $S' \subset S$ such that $\Delta = conv(S')$ is a simplex and $z$ is in the interior of $conv(S')$. We can assume that
$conv(S')$ is a minimal (by inclusion) simplex such that $S' \subset S$ and $z$ is in the interior of $conv(S')$.
 Then the interior of $\Delta $ does not contain points from $S$. 
Clearly, there is $s \in S$ such that $S' = S \setminus s$. Let $L$ be the line containing $s$ and $z$ and 
$sLt$ be the closed interval in $L$ such that $z$ is in the interior of $sLt$ and $t$ belongs to a face of 
$\Delta $. Let $sLt'$ be the maximal closed interval in $L$ that has no interior point of $\Delta $.
Obviously, there exist faces  $F$ and $F'$ of  $\Delta $ containing $t$ and $t'$, respectively. 
In plane language, $t'$ and $t$ are the first and the last common points of the ray $sL$ with $\Delta $.
%
%
Since $z$ is in a general position with respect to $S$, clearly $t$ and $t'$ are the interior points of $F$ and $F'$, respectively, 
and $dim (R(F)) = d-1$, 
and so $v(F)  = d$.
%
%
Let $s' = S' \setminus V(F)$. Then 
$\Delta ' = conv(F \cup s) = conv (S \setminus s')$ is a $z$-containing simplex, and so 
$V(F)  \cup s$ is a minimal $z$-containing subset of $S$.
By the above definition, 
 $\Delta =  conv (S \setminus s) = conv (F \cup s')$ is another $z$-containing simplex, and so 
 $V(F) \cup s' = S \setminus s$ is another minimal $z$-containing subset of $S$.
 
We claim that  $S \setminus x$ is a  $z$-avoiding subset of $S$ for every $x \in V(F)$. 
To prove this claim we consider two cases.
Suppose first that $x \in V(F) \setminus V(F')$. Then 
$S \setminus x = V(F') \cup s$ and by the definition of $F'$,  $z$ is not in the interior of $conv(V(F') \cup s)$. Therefore $S \setminus x$ is $z$-avoiding.
Now suppose  that $x \in V(F) \cap V(F')$. Then 
$S \setminus x = V(F'') \cup s$, where $F''$ is a face of $\Delta $ distinct from $F$ and $F'$. Obviously, $L \cap F'' = \emptyset $. Therefore $z$ is not an interior point of $ Conv(V(F'') \cup s)$, and so again $S \setminus x$ is $z$-avoiding.
 
Obviously, if $S \setminus x$ is  $z$-avoiding, then  $S \setminus x$ is also maximal $z$-avoiding. Therefore 
$S \setminus x$ is a maximal $z$-avoiding subset of $S$ for every $x \in V(F)$.
%
%
Also $V(F)$ is a $z$-avoiding subset of $S$ and since both $V(F) \cup s$ and $V(F) \cup s'$ are $z$-containing, clearly $V(F)$ is a maximal $z$-avoiding subset of $S$.
Thus, 
${\cal A}(S) = \{S \setminus x: x \in V(F)\} \cup \{V(F)\}$.
Also both $\Delta = conv(F \cup s')$ and 
$\Delta ' = conv(F \cup s)$ are  $z$-containing simplexes, and  so  ${\cal C}(S) = \{S \setminus s, S \setminus s'\}$.
Therefore, 
 $d +1= |{\cal A}(S)| = d~|{\cal C}(S)| - d +1 = d+1$.
 It follows that $(a3)$  holds.
 \\[1.5ex]
\indent
{\bf (p4)} The proof of  $(a4)$ is similar to that in ${\bf (p3)}$ because it can be checked that the inequality holds if
$|S| = d +3$. Claim $(a4)$ is also a particular case of Theorem \ref{mm} below.
\ep
\\[2ex]
\indent
Our next goal is to prove a strengthening of Theorem \ref{m} that takes into consideration the size of $S$.
\begin{lemma}
\label{good-vertex-in-conv(S)}
Let $S$ be a finite set of points in  $\mathbb{R}^d$, 
and $z$  in a general position with respect to $S$
{\em (and so $z \in \mathbb{R}^d \setminus S$)}. Let $P = conv (S)$ and $V = V(P)$.
Suppose that $|S| \ge d+2$ and $S$ is a $z$-containing set.
Then there exist a simplex $\Delta $ and 
$u \in S$
such that
\\[0.7ex]
$(a1)$ $V(\Delta ) \subseteq V$,
\\[0.7ex]
$(a2)$ $z$ is an interior point of $\Delta $,
\\[0.7ex]
$(a3)$ $\Delta $ and $P$  have a common facet hyperplane $H$, and 
\\[0.7ex]
$(a4)$ $u \in H$ and 
$S \setminus u$ is a $z$-containing set.
\end{lemma}
\bp (uses Lemma \ref{SimplexInPcontainingZero}).
We will consider two cases: $P$ is a simplex or 
$P$ is not a simplex.
\\[1.5ex]
\indent
${\bf (p1)}$ 
Suppose  that $P$ is a simplex.
Let $\Delta = P$.
Since $|S| \ge d+2$, there exists 
$v \in S \setminus V(P)$.

Suppose that $v$ is in the interior of $P$. 
Then obviously, there exists a facet $F$ of $P$ such that $z$ is in the interior of $conv(F \cup s)$. 
Let $u$ be a (unique) vertex of $V(P)$ not belonging to 
$F$.  Then $(\Delta , u)$ satisfies  $(a1)$ - $(a4)$, where $H$ is the  hyperplane containing $F$.

Now suppose that $v$ is not  in the interior of 
$P$. Then there exists a facet $T$ of $P$ containing $v$.  Put $u = v$. Then again $(\Delta , u)$ satisfies  $(a1)$ - $(a4)$, where $H$ is the  hyperplane containing $T$.
\\[1.5ex]
\indent
${\bf (p2)}$ Finally, suppose that $P$ is not a simplex.
By  Lemmas \ref{SimplexInPcontainingZero}, there exists a simplex $\Delta _z$ such that 
$V(\Delta _z) \subseteq V$ and $z$ is in the interior of 
$\Delta _z$.
Since $P$ is not a simplex, $\Delta _z \ne P$.
Therefore  $\Delta _z$ has a facet $F_z$ which is not a facet of $P$. Let $v$ be a unique vertex of  $\Delta _z$ such that $v \not  \in F_z$.
Let $L$ be the line containing $v$ and $z$. Let $vLz'$ be a maximal segment in $L$ such that $z \in vLz'$ and $vLz' \subset P$. Obviously, such segment exists (and is unique) and $z'$ belongs to a face $F$ of $P$.
Moreover, since  $z$ is in a general position with respect to $S$, clearly $F$ is a facet of $P$ and point $z'$  is in a general position with respect to $S' = V(F)$  (and so $S'$ is a $z'$-containing set) in $\mathbb{R}^{d-1}$.
By Lemma \ref{SimplexInPcontainingZero}, $F$ contains an $(d-1)$-dimensional simplex 
$\Delta '$ such that $z'$ is in the interior of $\Delta '$.
Then $\Delta = conv(\Delta ' \cup v)$ is a $d$-dimensional simplex satisfying $(a1)$, $(a2)$, and $(a3)$ with 
$H = R(F) = R(\Delta ')$.
Since $F$ is a facet of $P$, clearly  $V(F) \subset V(P) \subseteq S$ and
$V(F) \setminus V(F_z) \ne \emptyset $.
Let $u$ be an arbitrary point of $V(F) \setminus V(F_z)$. Obviously, $u \ne v$, for otherwise, 
$z \in vLz' \subseteq F$, and so  $z \in F \subset H$. This is impossible because $z$ is in the interior of $P$. 
Thus, $\Delta \subset  conv(S \setminus u)$, and so 
$S \setminus u$ is $z$-containing. 
\ep
\\[2ex]
\indent
Now we are ready to prove the following strengthening of Theorem \ref{m}. 
\begin{theorem}
\label{mm}
Let $S$ be a finite set of points in the 
$d$-dimensional space $\mathbb{R}^d$, $z \in \mathbb{R}^d \setminus S$.
Suppose that $S$ is $z$-containing, $z$  is in a general position with respect to $S$, and 
$|S| \ge d + 2$.
Then 
$|{\cal A}(S)| \le d~|{\cal C}(S)| - |S| + 3$.
\end{theorem}

\bp (uses Lemmas \ref{lemma!}, 
\ref{good-vertex-in-conv(S)}, and Theorem \ref{m}).
We prove our claim by induction on $|S|$.
If $|S| = d + 2$, then by Theorem \ref{m} $(a3)$, the claim is true. We assume  that our claim is true for $|S| = k \ge d+2$ and  prove that it is also true if $|S| = k +1$.
By Lemma \ref{lemma!},
\\[1ex]
\indent
 $|{\cal A}(S)| - |{\cal A}(S \setminus s)| = |{\cal A}_s(S)| \le 
 |{\cal H}^e_s(S)| \le |{\cal H}_s(S)| \le  
 d~|{\cal S}mpl_s(S)| = d~|{\cal C}_s(S)|$
 \\[1ex]
 for every $s \in S$.
 By Lemma \ref{good-vertex-in-conv(S)}, 
 there exist a point $u$ in $S$ and a simplex $\Delta $ in ${\cal S}impl_u(S)$ such that 
 $S \setminus u$ is $z$-containing  and 
$R(T) = R(F)$ for some facets $T$ and $F$ of $\Delta $ and $conv(S)$, respectively. Then obviously, $R(T)$ is a non-essential hyperplane in ${\cal H}_u(S)$, and therefore 
\\[0.5ex] 
$|{\cal H}^e_u(S)| \le |{\cal H}_u(S)| - 1$. 
Therefore 
\\[1ex]
\indent
$|{\cal A}(S)| - |{\cal A}(S \setminus u)| = |{\cal A}_u(S)| \le 
 |{\cal H}^e_u(S)| < |{\cal H}_u(S)| \le  d~|{\cal C}_u(S)|$.
\\[1ex]
By the induction hypothesis, we have:
$|{\cal A}(S \setminus u)| \le 
d~|{\cal C}(S\setminus u)| - |S \setminus u| + 3$.
\\[1ex]
Thus, our inductive step follows from the last two inequalities.
\ep

\end{document}